\documentclass[11pt]{amsart}
\usepackage{amssymb, latexsym, amsmath, amsfonts, tikz}
\usepackage{hyperref, enumitem, fancyhdr}

\newtheorem{thm}{Theorem}[section]

\theoremstyle{definition}

\theoremstyle{remark}

\numberwithin{equation}{section}
\theoremstyle{remark}

\newcommand{\mbb}{\mathbb}
\newcommand{\ra}{\rightarrow}

\newcommand{\pa}{\partial}
\newcommand{\ov}{\overline}
\newcommand{\sm}{\setminus}

\newcommand{\no}{\noindent}
\newcommand{\al}{\alpha}

\newcommand{\ti}{\tilde}

\newcommand{\be}{\beta}

\begin{document}
\title{On Grauert's examples of complete K\"{a}hler metrics}
\keywords{Holomorphic sectional curvature, complete K\"{a}hler metric}
\subjclass{Primary: 32Q05, 32Q10  ; Secondary : 32Q02}
\author{Sahil Gehlawat and Kaushal Verma}

\address{SG: Department of Mathematics, Indian Institute of Science, Bangalore 560 012, India}
\email{sahilg@iisc.ac.in}

\address{KV: Department of Mathematics, Indian Institute of Science, Bangalore 560 012, India}
\email{kverma@iisc.ac.in}

\begin{abstract}
Grauert showed that the existence of a complete K\"{a}hler metric does not characterize domains of holomorphy by constructing such metrics on the complements of complex analytic sets in a domain of holomorphy. In this note, we study the holomorphic sectional curvatures of such metrics in two prototype cases namely, $\mbb C^n \sm \{0\}, n \ge 2$ and $\mbb B^N \sm A$, $N \ge 2$ and $A \subset \mbb B^N$ is a hyperplane of codimension at least two. This is done by computing the Gaussian curvature of its restriction to the leaves of a suitable holomorphic foliation of these two examples. We also examine this metric on the punctured plane $\mbb C^{\ast}$ and show that it behaves very differently in this case.
\end{abstract}

\maketitle

\section{Introduction}

\no In addition to showing that every domain of holomorphy admits a real analytic, complete K\"{a}hler metric, Grauert \cite{Gr} showed that the converse is false in general by constructing such metrics on the complement of a closed complex analytic set in a given domain of holomorphy. This was done by first explicitly producing such a metric on 
$(\mbb C^{n})^{\ast} = \mbb C^n \sm \{0\}, n \ge 1$ and then pulling it back to a domain of the form $D \sm A$ (where $D \subset \mbb C^n$ is a domain of holomorphy and $A \subset D$ is a complex analytic set) by means of the defining equations for $A$. The possible degeneracies that may occur in pulling it back are taken care of by adding a suitable metric to it.

\medskip

The primary purpose of this note is to study the holomorphic sectional curvature of this metric on $(\mbb C^{n})^{\ast}$ where $n \ge 2$. In what follows, $n \ge 2$ unless stated otherwise. Furthermore, as an example, we apply Grauert's construction to the pair $D = \mbb B^N$, the unit ball in $\mbb C^N$, $N \ge 2$ and $A = \{z_1 = z_2 = \ldots = z_n =  0\}$, $2 \le n \le N$, to get such a metric on $\mbb B^N \sm A$, and also obtain information about its holomorphic sectional curvature. 

\medskip

To recall the construction of this metric on $(\mbb C^{n})^{\ast}$ (in addition to \cite{Gr}, see \cite{JP} for instance), note that the real analytic function $u : (0, \infty) \ra \mbb R$ defined by $u(t) = (t-1)/(t \log t)$ is positive on $(0, \infty)$ and satisfies
\[
\int_0^1 t u^2(t) \; dt < \infty \; \text{and} \; \int_0^1 u(t) \; dt = \infty.
\]
as a direct calculation shows. For $t \in (0, \infty)$, set $v(t) = \int_0^{t} \tau u^2(\tau) \; d \tau$ and $U(t) = \int_0^t v(\tau)/\tau \; d \tau$ so that $U'(t) = v(t)/t$. Finally, define $H(z) = U(\vert z \vert^2)$ and let $g_{\al \ov \be} =  \pa^2 H/\pa z_{\al} \pa \ov z_{\be}$ for $1 \le \al, \be \le n$. A computation shows that
\[
g_{\al \ov \be}(z) = U''(\vert z \vert^2) \ov z_{\al} z_{\be} + U'(\vert z \vert^2) \delta_{\al \be}
\]
and hence for a tangent vector $X = (X_1, X_2, \ldots, X_n) \in \mbb C^n$ based at $z \in (\mbb C^{n})^{\ast}$,
\[
\sum_{\al, \be} g_{\al \ov \be}(z) X_{\al} \ov X_{\be} = \sum_{\al, \be} \left(  U''(\vert z \vert^2) \ov z_{\al}  z_{\be} + U'(\vert z \vert^2) \delta_{\al \be} \right) X_{\al} \ov X_{\be}.
\]                                     
The right side simplifies as
\[
\left( u^2(\vert z \vert^2) - \frac{v(\vert z \vert^2)}{\vert z \vert^4}  \right) | \langle X, z \rangle |^2 + \frac{v(\vert z \vert^2)}{\vert z \vert^2} \vert X \vert^2 = u^2(\vert z \vert^2) |\langle X, z \rangle |^2 + \frac{v(\vert z \vert^2)}{\vert z \vert^2} \left( \vert X \vert^2 - \frac{| \langle X, z \rangle |^2}{\vert z \vert^2} \right)
\]
where $\langle X, z \rangle =  X_1 \ov z_1 + X_2 \ov z_2 + \ldots + X_n \ov z_n$ is the standard Hermitian inner product on $\mbb C^n$. The Cauchy-Schwarz inequality then shows that the prescription
\[
\tilde g(z, X) = \sum_{\al, \be} g_{\al \ov \be}(z) X_{\al} \ov X_{\be} \ge 0
\]
and hence defines a real analytic K\"{a}hler pseudometric, with global potential $H$, on $(\mbb C^{n})^{\ast}$. By replacing $g_{\al \ov \be}(z)$ with $g_{\al \ov \be}(z) + \delta_{\al \be}$ amounts to working with $g(z, X) = \tilde g(z, X) + \vert X \vert^2$ which is clearly positive and hence produces the desired metric. We will refer to $g$ as the Grauert metric on $(\mbb C^{n})^{\ast}$ and record the expression
\[
g(z, X) = \left( 1 + \frac{v(\vert z \vert^2)}{\vert z \vert^2} \right) \vert X \vert^2 + \left( u^2(\vert z \vert^2) - \frac{v(\vert z \vert^2)}{\vert z \vert^4} \right) \vert \langle X, z \rangle \vert^2
\]
for future use. Note that $g(z, X)$ is rotation invariant in the sense that $g(z, X) = g(\mbb U z, \mbb U X)$ for every unitary rotation $\mbb U$.

\medskip

By \cite{Wu}, recall that the holomorphic sectional curvature of $g$ at $z \in (\mbb C^{n})^{\ast}$ along a tangent direction $X$ can be recovered as the supremum of the Gaussian curvatures of all possible nonsingular holomorphic disks through $z$ which are tangent to $X$. In other words, let $\phi$ be a holomorphic map from a neighbourhood of the origin in $\mbb C$ into 
$(\mbb C^{n})^{\ast}$ with $\phi(0) = z$ and $d\phi(0) = X$. Let $K_{\phi^{\ast}g}(0)$ be the Gaussian curvature of the pull back metric $\phi^{\ast}g$ which is well defined near the origin. Then, $K_g(z, X)$, the holomorphic sectional curvature of $g$ at $z$ along $X$ satisfies
\[
K_g(z, X) = \sup K_{\phi^{\ast}g}(0)
\]
where the supremum is taken over all such $\phi$. Let $K^{\pm}_g(z)$  be the supremum and infimum respectively of $K_g(z, X)$ over all tangent vectors $X$ at $z$. 

\medskip

Let $\phi$ be such an embedding from a neighbourhood of the origin into $(\mbb C^{n})^{\ast}$ with $\phi(0) = z$ and $d\phi(0) = X$. For a given unitary rotation $\mbb U$, the map $\mbb U \circ \phi$ defines an embedding at $\mbb U \circ \phi(0) = \mbb U z$ with $d(\mbb U \circ \phi)(0) = \mbb U X$. Note that
\[
(\mbb U \circ \phi)^{\ast} g = (\phi^{\ast} \circ \mbb U^{\ast}) g = \phi^{\ast} g
\]
since $g$ is rotation invariant. Therefore, 
\[
K_{\phi^{\ast} g}(0) = K_{(\mbb U \circ \phi)^{\ast} g}(0)
\]
which shows that $K_g(z, X) = K_g(\mbb U z, \mbb U X)$ and hence $K^{\pm}_g(z) = K^{\pm}_g(\mbb U z)$. Since this holds for every $\mbb U$, $K^{\pm}_g$ are constant on spheres centered at the origin. 

\medskip

As $(\mbb C^{n})^{\ast}$, $n \ge 2$,  is not a domain of holomorphy, Shiffman's theorem \cite{Sh} implies that $K^+_g$ cannot be non-positive everywhere. Thus, there must be points arbitrarily close to the origin and holomorphic tangent vectors $X$ based at them along which $K_g(X) > 0$. This was strengthened by Balas \cite{Ba} who showed that {\it every} rotation invariant, complete K\"{a}hler metric on $(\mbb C^{n})^{\ast}$ admits holomorphic tangent vectors, based at points arbitrarily close to the origin, along which it is positively curved. This was done by considering tangent vectors that are orthogonal to the radial direction. However, it turns out that much more can be said for the Grauert metric $g$ by computing the curvature of its restriction to the leaves of the holomorphic vector field 
\begin{equation}
X = X(z) = z_1 \frac{\pa}{\pa z_1} + z_2 \frac{\pa}{\pa z_2} + \ldots + z_{n-1} \frac{\pa}{\pa z_{n-1}} + \al z_n \frac{\pa}{\pa z_n}, \;\; \al < 0
\end{equation}
which has an isolated singularity at the origin. Let $\mathcal F^n_X$ be the foliation induced by $X$. For $z \in (\mbb C^n)^{\ast}$, let $\mathcal L_z \subset \mathcal F^n_X$ be the leaf passing through $z$. Let $\kappa(z)$ be the Gaussian curvature of $g|\mathcal L_z$, the restriction of $g$ to $\mathcal L_z$, at $z$. By estimating $\kappa$ near the origin and noting that $K^+_g \ge \kappa$ provides a lower bound for $K^+_g$ near the origin, and a similar calculation, that is independent of this specific $X$, yields information on $K_g(z, X)$ and hence $K^-_g$ far away from the origin.

\begin{thm}
For the Grauert metric $g$ on $(\mbb C^{n})^{\ast}$, where $n \ge 2$, $K^+_g(z) \ra +\infty$ as $\vert z \vert \ra 0$. Furthermore, there exists $R > 0$ such that $K^-_g(z) < 0$ for all
$\vert z \vert \ge R$.
\end{thm}

As alluded to earlier, we will now consider a specific example of a complete metric that arises by applying Grauert's construction to $\mbb B^N \sm A$ where $A = \{z_1 = z_2 = \ldots = z_n =  0\}$, $2 \le n \le N$ and $N \ge 2$. The Grauert metric $g$ on $(\mbb C^n)^{\ast}$ can be pulled back via the natural projection $\pi : \mbb B^N \sm A \ra (\mbb C^n)^{\ast}$
\[
\pi(z_1, z_2, \ldots, z_N) = (z_1, z_2, \ldots, z_n)
\] 
to get $\pi^{\ast} g$ which is a K\"{a}hler pseudometric on $\mbb B^N \sm A$. By adding to this a complete metric on $\mbb B^N$ say, $\tilde g$, gives $\pi^{\ast} g + \ti g$ which is complete on $\mbb B^N \sm A$. Note that the choice of $\ti g$ is not unique and as a prototype, we will consider 
the Bergman metric
\[
\ti g_B = (N+1) \left( \frac{\vert \langle z, dz \rangle \vert^2}{(1 - \vert z \vert^2)^2} + \frac{\vert dz \vert^2}{1 - \vert z \vert^2} \right)
\]
which is a complete K\"{a}hler, constant curvature example. Thus, $g_B = \pi^{\ast} g + \ti g_B$
is an example of a complete K\"{a}hler metric on $\mbb B^N \sm A$ and the corresponding functions $K^{\pm}_{g_B}$ are well defined. Note that $X$, when viewed as a vector field on 
$\mbb B^N \sm A$, is singular along $A$. Let $\mathcal F^N_X$ be the foliation induced by $X$ on $\mbb B^N \sm A$. Fix $w = (0, \ldots, 0, w_{n+1}, \ldots, w_N) \in A$, let $r^2 = 1- \vert w \vert^2$ and consider the slice
\[
D_w = \{ z \in \mbb B^N : z_j = w_j, \; n+1 \le j \le N \}.
\]
Note that $D_w \sm A$ can be identified with a punctured ball in $\mbb C^n$ of radius $r$. Also, for $z \in D_w$, the entire leaf $\mathcal L_z \subset \mathcal F^N_X$ passing through $z$ is contained in $D_w$. Let $\mathcal K(z)$ be the Gaussian curvature of $g_B | \mathcal L_z$, the restriction of $g_B$ to $\mathcal L_z$, at $z$. Again, as $K^+_{g_B}(z) \ge \mathcal K(z)$ , it suffices to estimate $\mathcal K(z)$ near $w$ in order to obtain a lower bound for $K^+_{g_B}(z)$ near $w$. In addition, $\mathcal K(z)$ can be estimated directly near the boundary of $D_w \sm A$.

\begin{thm}
For the metric $g_B$ on $\mbb B^N \sm A$, $K^+_{g_B}(z) \ra +\infty$ as $z \ra w$ along $D_w$. Furthermore, $\mathcal K(z) \ra -6/(N+1)$ as $\vert z \vert \ra 1$ along $D_w$.
\end{thm}

The concluding section contains a discussion of this metric on $\mbb C^{\ast}$, i.e., when $n = 1$, as well as miscellaneous remarks related to the computations that follow. Other aspects of complete metrics on complements of analytic sets can be found in \cite{DF1}, \cite{DF2} and \cite{Oh}.


\section{Proof of Theorem 1.1}

\no We will begin by computing a general expression for the Gaussian curvature of the leaves of a foliation induced by an arbitrary holomorphic vector field $X$ on $\mbb C^n$ with $\text{sing}(X) = \{0\}$. For 
$z \in  (\mbb C^{n})^{\ast}$, let $Z = Z(T)$ be the solution of
\[
\dot{Z} = dZ/dT = X(Z(T)), \; Z(0) = z
\]
where $Z$ is defined in a neighbourhood of the origin in the plane. Since $\dot{Z} = X(Z)$, the pull back metric $Z^{\ast}g(T) = h(T) \vert dT \vert^2$ where
\begin{equation}
h(T) = \left( 1 + \frac{v(\vert Z \vert^2)}{\vert Z \vert^2} \right) \vert X(Z) \vert^2 + \left( u^2(\vert Z \vert^2) - \frac{v(\vert Z \vert^2)}{\vert Z \vert^4} \right) \vert \langle Z, X(Z) \rangle \vert^2.
\end{equation}
Its curvature
\begin{equation}
\kappa(z) = -2 (h(T))^{-1}\pa \ov \pa \log h(T) \big|_{T = 0} = -2 (h(0))^{-3} \left( h(0) \; \pa \ov \pa h(0) - \pa h(0) \; \ov \pa h(0) \right)
\end{equation}
which will be simplified by computing expressions for $\pa h, \ov \pa h$ and $\pa \ov \pa h$. For this, abbreviate $u(\vert Z \vert^2), v(\vert Z \vert^2)$ by just $u, v$ and likewise for their derivatives with the understanding that they are functions of $\vert Z \vert^2$. Also, note that
\begin{equation}
\pa \vert Z \vert^2 = \langle X(Z), Z \rangle, \;\; \ov \pa \vert Z \vert^2 = \langle Z, X(Z) \rangle, \;\; \pa X(Z) = DX(Z)(X(Z)) \; \text{and} \; \ov \pa X(Z) = 0.
\end{equation}
Repeated use of these identities in simplifying the resulting expressions shows that
\[
\pa h = f_1(T) + f_2(T) + f_3(T) \;\; \text{and} \;\; \ov \pa h = \ov{f_1(T)} + \ov{f_2(T)} + \ov{f_3(T)}
\]
where
\begin{align*}
f_1(T) &= \left( u^2 - \frac{v}{\vert Z \vert^4} \right) \left( 2 \vert X(Z) \vert^2 \langle X(Z), Z \rangle + \langle DX(Z)(X(Z)), Z \rangle \langle Z, X(Z) \rangle \right), \\
f_2(T) &= \left( 1 + \frac{v}{\vert Z \vert^2} \right) \langle DX(Z)(X(Z)), X(Z) \rangle, \; \text{and} \\
f_3(T) &= \left( 2uu' - \frac{u^2}{\vert Z \vert^2} + \frac{2v}{\vert Z \vert^6}  \right) \langle X(Z), Z \rangle \vert \langle X(Z), Z \rangle \vert^2.
\end{align*}
Differentiating once again 
\[
\ov \pa f_1(T) = \left( 2uu' -\frac{u^2}{\vert Z \vert^2} + \frac{2v}{\vert Z \vert^6} \right) A_1  + \left(u^2 - \frac{v}{\vert Z \vert^4} \right)  \left(A_2 + 2 \vert X(Z) \vert^4  + \vert \langle DX(Z)(X(Z)), Z \rangle \vert^2 \right)
\]
where the $A_i$'s are multi-term expressions that can be explicitly written, but their exact form is not important. The relevant point is that $A_1, A_2$ both contain either $\langle X(Z), Z \rangle$ or its conjugate as a factor. Similarly,
\[
\ov \pa f_2(T) = \left( u^2 - \frac{v}{\vert Z \vert^4} \right) A_3 + \left( 1 + \frac{v}{\vert Z \vert^2}  \right) \vert DX(Z)(X(Z)) \vert^2
\]
where $A_3$ contains $\langle Z, X(Z) \rangle$ as factor and finally, it can be seen that every term in $\ov \pa f_3(T)$ also contains either $\langle X(Z), Z \rangle$ or its conjugate as a factor. Thus, $\ov \pa \pa h = \ov \pa f_1(T) + \ov \pa f_2(T) + \ov \pa f_3(T)$ can be expressed as the sum of several terms that can be grouped depending on whether they contain $\langle X(Z), Z \rangle$ (or its conjugate) or not.

\medskip

Let $\mathcal S = \{ z \in \mbb C^n \sm \{0\} : \langle z, X(z) \rangle = 0 \}$. This is precisely the set where the vector field $X$ is orthogonal to the radial direction. Note that $\mathcal S$ may be empty as the example of $X = (z_1, z_2, \ldots, z_n)$ shows. But assuming that $\mathcal S \not= \emptyset$, note that for $Z(0) = z  \in \mathcal S$,  
\[
h(0) = \left( 1 + \frac{v}{\vert z \vert^2} \right) \vert X(z) \vert^2, \; f_2(0) = \left( 1 + \frac{v}{\vert z \vert^2} \right) \langle DX(z)(X(z)), X(z)  \rangle
\]
and $f_1(0) = f_3(0) = 0$. Hence
\[
\pa h(0) = \left( 1 + \frac{v}{\vert z \vert^2}  \right) \langle DX(z)(X(z)), X(z) \rangle = \ov{\ov \pa h(0)}
\]
Further, on $\mathcal S$, all the $A_i$'s that occur in the expressions for $\ov \pa f_i$ vanish and hence
\[
\ov \pa f_1(0) = \left( u^2 - \frac{v}{\vert z \vert^4}  \right) \left( 2 \vert X(z) \vert^4 + \vert \langle DX(z)(X(z)), z \rangle \vert^2 \right),  \; \; \ov \pa f_2(0) = \left( 1 + \frac{v}{\vert z \vert^2} \right) \vert DX(z)(X(z)) \vert^2
\]
and $\ov \pa f_3(0) = 0$, and thus
\[
\ov \pa \pa h(0) = \left( u^2 - \frac{v}{\vert z \vert^4}  \right) \left( 2 \vert X(z) \vert^4 + \vert \langle DX(z)(X(z)), z \rangle \vert^2 \right) +  \left( 1 + \frac{v}{\vert z \vert^2} \right) \vert DX(z)(X(z)) \vert^2.
\]
Putting all this together in (2.2) and writing $\vert V \wedge W \vert^2 = \vert V \vert^2 \vert W \vert^2 - \vert \langle V, W \rangle \vert^2$ for vectors $V, W \in \mbb C^n$ leads to $\kappa(z) = N(z) / D(z)$ where $N(z) = -2(A(z) + B(z))$ with
\begin{align}
A(z) &= (\vert z \vert^4 u^2 - v)( 2\vert X(z) \vert^6 + \vert \langle DX(z)(X(z)), z \rangle \vert^2 \vert X(z) \vert^2 ), \\
B(z) &= (\vert z \vert^2 + v) \vert z \vert^2 \vert X(z) \wedge DX(z)(X(z)) \vert^2  
\end{align}
and 
\begin{equation}
D(z) = (\vert z \vert^2 + v)^2 \vert X(z) \vert^6
\end{equation}
after some simplification.

\medskip

Note that $B(z) \ge 0$. Further, the sign of $A(z)$, and hence $\kappa(z)$, is determined by looking at 
$(\vert z \vert^4 u^2 - v)$. Setting $t = \vert z \vert^2$, the derivative of the function $\eta : (0, \infty) \ra \mbb R$ defined by
\[
\eta(t) = t^2u^2(t) - v(t)
\]
satisfies 
\[
\lim_{t \ra \infty} \eta'(t) = \lim_{t \ra \infty} \left( t u^2(t)  + 2t^2u(t) u'(t) \right) = \lim_{t \ra \infty} \left ( \frac{(t-1)^2 \log t}{t(\log t)^3} \left( 1 + 2 \left( \frac{1}{t-1} - \frac{1}{\log t} \right) \right) \right) = +\infty.
\]
Thus, $\eta(t) \ra +\infty$ as $t \ra \infty$ and hence there exists $R > 0$ such that $A(z) > 0$, and consequently $\kappa(z) < 0$, for
$z \in \mathcal S \cap \{ \vert z \vert^2 > R^2 \}$. 


\subsection{Behaviour near infinity:} Two observations can be made about the calculations so far. One, the critical value of $R$ beyond which $\kappa < 0$ does {\it not} depend on the vector field $X$ as along as $z \in \mathcal S$ -- it only depends on the functions $u, v$ that define the Grauert metric $g$, and second, all these calculations are purely local, i.e., what is needed is a non-singular holomorphic vector field that is defined in some neighbourhood contained in $(\mbb C^n)^{\ast}$.

\medskip

Pick a $z_0$ with $\vert z_0 \vert > R$ and let $V$ be a holomorphic tangent vector at $z_0$. By Wu \cite{Wu}, there are local holomorphic coordinates $(w_1, w_2, \ldots, w_n)$ defined in a neighbourhood $\Omega$ of $z_0 = 0$ in which $V = \pa / \pa w_n |_{0}$ and the Gaussian curvature of the restriction of the Grauert metric $g$ to the holomorphic disc $\{w_1 = w_2 = \ldots = w_{n-1} = 0\}$ equals $K_g(z_0, V)$. The assignment $p \mapsto \pa / \pa w_n |_p$ defines a non-vanishing holomorphic vector field $X(p)$ in $\Omega$ and the leaves of the induced foliation are given by $\{w_1 = c_1, w_2 = c_2, \ldots, w_{n-1} = c_{n-1} \}$ for small values of $c_i \in \mbb C$. Note that $X(0) = V$. Now suppose, in addition, that $\langle z_0, V \rangle = 0$. The calculations above show that $K_g(z_0, V) = \kappa(z_0) < 0$ and hence $K^-_g(z_0) \le K_g(z_0, V) < 0$ and by rotation invariance, it follows that $K^-_g(z) < 0$ for all $\vert z \vert = \vert z_0 \vert$.


\subsection{Behaviour near the origin:} To estimate $K^+_g$ near $z = 0$, the vector field $X$ in (1.1) will be used. Write $X(z) = (z_1, z_2, \ldots, \alpha z_n)$ where $\al < 0$ and compute $DX(z)(X(z)) = (z_1, z_2, \ldots, \alpha^2 z_n)$. For $n \ge 3$,  
\[
\mathcal S = \{ z \in \mbb C^n \sm \{0\} : \langle z, X(z) \rangle = \vert z_1 \vert^2 + \vert z_2 \vert^2 + \ldots \vert z_{n-1} \vert^2 + \alpha \vert z_n \vert^2 = 0\}
\]
is a smooth quadric whose Levi form has exactly one zero eigenvalue and $(n-2)$ non-zero eigenvalues everywhere.  When $n = 2$,
\[
\mathcal S = \{ z \in \mbb C^2 \sm \{0\} : \langle z, X(z) \rangle =\vert z_1 \vert^2 + \alpha \vert z_2 \vert^2 = 0 \}
\]
 and is smooth Levi flat, but in either case, the origin $0 \in \ov {\mathcal S}$ and $\mathcal S \cup \{0\}$ is a real algebraic hypersurface with an isolated singularity at the origin. On 
$\mathcal S$, the various terms in (2.4) - (2.6) simplify as
\[
\vert z \vert^2 = (1- \al) \vert z_n \vert^2, \; \; \vert X(z) \vert^2 = (\al^2 - \al) \vert z_n \vert^2, \;\; \vert DX(z)(X(z) \vert^2 = (\al^4 - \al) \vert z_n \vert^2
\]
and 
\[
\langle DX(z)(X(z)), z \rangle = (\al^2 - \al) \vert z_n \vert^2, \;\; \langle DX(z)(X(z)), X(z) \rangle = (\al^3 - \al) \vert z_n \vert^2.
\]
Further, 
\begin{align*}
\vert z \vert^2 \vert X(z) \wedge DX(z)(X(z)) \vert^2 &= (1 - \al) \vert z_n \vert^2 \left( (\al^2 - \al)(\al^4 - \al) \vert z_n \vert^4 - (\al^3 - \al)^2 \vert z_n \vert^4  \right) \\
                                                                                       &= (\al^2 - \al)^3 \vert z_n \vert^6.
\end{align*}
Therefore, $A(z), B(z)$ and $D(z)$ can be determined and hence
\begin{equation}
\kappa(z) = - \frac{2\left( 3 \vert z \vert^4 u^2(\vert z \vert^2) - 2 v(\vert z \vert^2) + \vert z \vert^2  \right)}{\left( \vert z \vert^2 + v(\vert z \vert^2) \right)^2}
\end{equation}
for $z \in \mathcal S$. The interesting thing here is that this expression is independent of $\al$. To analyze $\kappa(z)$, write $\vert z \vert^2 = t$ and consider $f : (0, \infty) \ra \mbb R$ defined by
\[
f(t) = - \frac{2(3t^2u^2(t) - 2 v(t) + t)}{(t + v(t))^2}
\]
which is clearly real analytic on $(0, \infty)$.

\medskip

\no {\it Claim:} $f(1) < 0$ and $\lim_{t \ra 0} f(t) = + \infty$.

\medskip

Recall that $u(t) = (t-1)/(t \log t)$ and $v(t) = \int_0^t \tau u^2(\tau) \; d \tau$. To show that $f(1) < 0$, it will suffice to prove that $v(1) < 2$ since $u(1) = 1$. Change variables by writing
 $\tau = e^{-x}$ to get
\[
v(1) = \int_0^1 \frac{(\tau - 1)^2}{\tau (\log \tau)^2} \; d \tau = \int_0^{\infty} \frac{(1 - e^{-x})^2}{x^2} \; dx = \int_0^1  \frac{(1 - e^{-x})^2}{x^2} \; dx + \int_1^{\infty}  \frac{(1 - e^{-x})^2}{x^2} \; dx
\] 
and note that $(1 - e^{-x})^2/x^2 \le 1$ on $[0, 1)$ and $(1 - e^{-x})^2 < 1$ on $[1, \infty)$ which respectively imply that the first term is at most $1$ and the second is strictly less than $1$. Thus, $v(1) < 2$.

\medskip

For the second claim, observe that both $u(t), tu^2(t) \ra +\infty$ and $t^2 u^2(t), v(t) \ra 0$ as $t \ra 0$. By substituting for $u, u', v$ and $v'$, and repeatedly using L'H\^{o}pital's rule, 
\[
\lim_{t \ra 0} f(t) = -\lim_{t \ra 0} \frac{2 \left( 4tu^2(t) + 6t^2 u(t) u'(t) + 1 \right)}{2t + 2tv(t)u^2(t) + 2v(t) + 2t^2 u^2(t)} = +\infty.
\]
Therefore, there exists a $c \in (0, 1)$ such that $f(c) = 0$. By the real analyticity of $f$ and the fact that $f(t) \ra +\infty$ near $t=0$, there can be only finitely many such $c$. Let $a > 0$ be the minimum of the zeros of $f$ in $(0,1)$. Hence, for $z \in \mathcal S \cap \{ \vert z \vert^2 < a \}$, $\kappa(z) > 0$ and this implies that $K^+_g(z) > 0$. By rotation invariance,
$K^+_g > 0$ everywhere in the punctured ball $\{ \vert z \vert^2 < a  \} \sm \{0\}$. In fact, $K^+_g(z) \ra +\infty$ as $\vert z \vert \ra 0$. The proof of Theorem 1.1 is complete.


\section{Proof of Theorem 1.2}

\no Recall that the metric $g_B = \pi^{\ast} g + \ti g_B$, where
\[
\pi^{\ast}g(z, X) = \left( 1 + \frac{v(\vert \pi(z) \vert^2)}{\vert \pi(z) \vert^2} \right) \vert \pi(X) \vert^2 + \left( u^2(\vert \pi(z) \vert^2) - \frac{v(\vert \pi(z) \vert^2)}{\vert \pi(z) \vert^4} \right) \vert \langle \pi(X), \pi(z) \rangle \vert^2 
\]
for $(z, X) \in \mbb (B^N \sm A) \times \mbb C^N$ and $\ti g_B$ is the Bergman metric on $\mbb B^N$.  For $z \in D_w$, let $Z = Z(T)$ be the solution of 
\[
\dot{Z} = dZ/dT = X(Z(T)), \; \; Z(0) = z \in D_w.
\]
Since $X$ as in (1.1) is being viewed as a vector field on $\mbb B^N \sm A$, it follows that $Z_j(T) \equiv w_j$ for $n+1 \le j \le N$. Also,
\[
\vert Z_1(T) \vert^2  + \vert Z_2(T) \vert^2 + \ldots + \vert Z_n(T) \vert^2 < 1 - \vert w \vert^2 = r^2.
\]
Furthermore, $\pi : D_w \sm A \ra B^n(0, r) \sm \{0\}$ is a biholomorphism and let
\[
\ti Z(T) = \pi(Z(T)) = (Z_1(T), Z_2(T), \ldots, Z_n(T)) \in \mbb B^n(0, r) \sm \{0\}
\]
and set $\ti z = \pi(z)$. For $\ti z = (z_1, z_2, \ldots, z_n) \in \mbb B^n(0, r)$, define
\[
\ti X(\ti z)  =  z_1 \frac{\pa}{\pa z_1} + z_2 \frac{\pa}{\pa z_2} + \ldots + z_{n-1} \frac{\pa}{\pa z_{n-1}} + \al z_n \frac{\pa}{\pa z_n}, \;\; \al < 0
\]
and denote the corresponding foliation by $\mathcal F_{\ti X}$. Then $\text{sing}(\ti X) = \{0\} = \pi(w)$ and note that $\ti X(\ti z) = \pi \circ X \circ \pi^{-1}(\ti z)$ and $\vert X(Z) \vert^2 = \vert \ti X(\ti Z) \vert^2$. Also, note that $\vert Z(T) \vert^2 = \vert \ti Z(T) \vert^2 + \vert w \vert^2$ which gives $1 - \vert Z(T) \vert^2 = r^2 - \vert \ti Z(T) \vert^2$ and $\langle Z(T), X(Z(T)) \rangle = \langle \ti Z(T), \ti X(\ti Z(T)) \rangle$.

\medskip

With $\mathcal S$ as before and $w \in A$, let 
\[
\mathcal S_w = \mathcal S \cap D_w = \{ z \in D_w \sm A : \vert z_1 \vert^2 + \vert z_2 \vert^2 + \ldots + \al \vert z_n \vert^2 = 0  \}
\]
which describes the locus along which the radial vector $z$ (along the slice $D_w$) is orthogonal to $\ti X$, i.e., $\langle \ti z, \ti X(\ti z) \rangle = \langle z, X(z) \rangle = 0$. As before, for $z \in \mathcal S_w$, 
\[
\vert \ti X(\ti z) \vert^2 \vert \langle D \ti X(\ti z)(\ti X(\ti z)), \ti z \rangle^2 = \vert \ti X(\ti z) \vert^6 = \vert \ti z \vert^2 \vert D \ti X(\ti z)(\ti X(\ti z)) \wedge \ti X(\ti z) \vert^2.
\]
Using the above observations, writing $\ti Z = \ti Z(T)$ and computing as before, $Z^{\ast} g_B(T) = h_B(T) \vert dT \vert^2$ where
\[
h_B(T) = \left( 1 + \frac{v(\vert \ti Z \vert^2)}{\vert \ti Z \vert^2} + \frac{N+1}{r^2 - \vert \ti Z \vert^2}\right) \vert \ti X (\ti Z) \vert^2 + \left( u^2(\vert \ti Z \vert^2) - \frac{v(\vert \ti Z \vert^2)}{\vert \ti Z \vert^4} + \frac{N+1}{(r^2 - \vert \ti Z \vert^2 )^2} \right) \vert \langle \ti Z, \ti X(\ti Z) \rangle \vert^2.
\]
The calculations leading to $\mathcal K(z)$, the curvature of $Z^{\ast} g_B(T)$, are similar in spirit to those in Theorem 1.1. Indeed, what is required are expressions for $h_B(0), \pa h_B(0), \ov \pa h_B(0)$ and $\pa \ov \pa h_B(0)$ when $Z(0) = z \in \mathcal S_w$. Of these, the easiest to write down is
\begin{equation}
h_B(0) = \left( 1 + \frac{v}{\vert \ti z \vert^2} + \frac{N+1}{r^2 - \vert \ti z \vert^2}  \right) \vert \ti X (\ti z) \vert^2
\end{equation}
where $u, v$ are understood to be functions of $\vert \ti z \vert^2$. For the others, several pages of messy calculations along similar lines as in Theorem 1.1 yield rather unwieldy expressions for $\pa h_B, \ov \pa h_B$ and $\pa \ov \pa h_B$ which simplify considerably when $Z(0) = z \in \mathcal S_w$ to give:
\begin{equation}
h_B(0) \; \pa \ov \pa h_B(0) - \pa h_B(0) \; \ov \pa h_B(0) = I(\ti z) \cdot J(\ti z)
\end{equation}
where
\begin{align}
I(\ti z) &= \left( \vert \ti z \vert^2 + v + \frac{(N+1) \vert \ti z \vert^2}{(r^2 - \vert \ti z \vert^2)}  \right) \frac{\vert \ti X(\ti z) \vert^6}{\vert \ti z \vert^6}, \; \text{and} \\
J(\ti z) &= \left( 3 \vert \ti z \vert^4 u^2 - 2v + \vert \ti z \vert^2 + \frac{(N+1)\vert \ti z \vert^2(r^2 + 2 \vert \ti z \vert^2)}{(r^2 - \vert \ti z \vert^2)^2} \right).
\end{align}
Finally, by combining (3.1) -- (3.4)
\[
\mathcal K(z) = - P(\ti z) / Q(\ti z)
\] 
where
\begin{align*}
P(\ti z) &= 2 \left( 3 \vert \ti z \vert^4 u^2 -2v + \vert \ti z \vert^2 + \frac{(N+1)\vert \ti z \vert^2(r^2 + 2 \vert \ti z \vert^2)}{(r^2 - \vert \ti z \vert^2)^2}  \right), \; \text{and}\\
Q(\ti z) &= \left( \vert \ti z \vert^2 + v + \frac{(N+1)\vert \ti z \vert^2}{(r^2- \vert \ti z \vert^2)}  \right)^2.
\end{align*}

\no To analyze $\mathcal K(z)$, write $\vert \ti z \vert^2 = t$ and define $f : (0, r^2) \ra \mbb R$ as
\[
f(t) = - \frac{2 \left( 3t^2 u^2(t) - 2v(t) + t + \frac{(N+1)t(r^2+2t)}{(r^2-t)^2}  \right)}{\left( t + v(t) + \frac{(N+1)t}{(r^2-t)}  \right)^2}
\]
which is real analytic on $(0, r^2)$. Theorem 1.2 is a consequence of the following:

\medskip

\no {\it Claim:} $f(r^2/2) < 0$, $\lim_{t \ra r^2}f(t) \ra -6/(N+1)$ and $\lim_{t \ra 0}f(t) \ra +\infty$.

\medskip

Note that the numerator in $f(r^2/2)$ is 
\[
-2 \left( 3(r^2/2)^2 u^2(r^2/2)  -2v^2(r^2/2) + r^2/2 + 4(N+1) \right)
\]
and this will be negative if it can be shown that $v(r^2/2) < 2(N+1)$. But since $r^2/2 < 1/2$, $v(1) < 2$ and $v$ is an increasing function on $(0, \infty)$, it follows that $v(r^2/2) < v(1) < 2 < 2(N+1)$.

\medskip

To determine the behaviour of $f$ near the origin, note that both the numerator and denominator in $f$ approach zero as $t \ra 0$ and repeated use of  L'H\^{o}pital's rule shows that $\lim_{t \ra 0}f(t) = +\infty$. Finally, note that near $t = r^2$, both $(3t^2u^2 - 2v + t)$ and $(t + v)$ are bounded while the remaining terms in $f$ approach $+\infty$. Simplifying them gives
\[
\lim_{t \ra r^2} f(t) = \lim_{t \ra r^2} -2(r^2 + 2t)/(t(N+1)) = -6/(N+1).
\]

When $n = N$, $A = \{0\}$ and $\mbb B^N \sm A$ is the punctured ball with the metric $g_B$. In this case, there is no need to work with slices for then $X$ as in (1.1) has an isolated singularity at the origin. All computations above remain valid in this case as well.

\section{Concluding Remarks}

\medskip

\no (i) When $n = 1$, note that
\begin{equation}
g(z, X) = \left( 1 + \vert z \vert^2 u^2(\vert z \vert^2) \right) \vert X \vert^2
\end{equation}
as a bit of simplification shows. The curvature properties of this metric on $\mbb C^{\ast}$ are very likely to be considered as part of folklore knowledge, but our inability to find any explicit reference only piqued our curiosity and made us include the following observations.

\medskip

The conformal density $h (z) = \left( 1 + \vert z \vert^2 u^2(\vert z \vert^2) \right) \ra +\infty$ as $\vert z \vert \ra 0$ and a calculation shows that 
the curvature of $g$ is $\kappa(z) = f(t), t = \vert z \vert^2$, where $f(t) = - P(t)/Q(t)$ with
\[
P(t) = 2 \left( u^2 + 2t^2(u')^2 + 2t^2 u u''+ 2t^3 u^3 u'' + 6tuu' + 2t^2u^3 u' - 2t^3 (uu')^2  \right)
\]
and $Q(t) = (1 + tu^2)^3$.

\medskip

\no {\it Claim:} $f \le 0$ on $(0, \infty)$, $\lim_{t \ra 0} f(t) = -4$ and $\lim_{t \ra \infty} f(t) = 0$. 

\medskip

\no To compute the limits of $f$ as $t \ra 0$ and $t \ra \infty$, write $P(t)  = 2(M_1(t) + M_2(t))$ where
\[
M_1(t) = 2t^3u^3u'' + 2t^2u^3u' - 2t^3u^2(u')^2 \; \text{and} \; M_2(t) =  u^2 + 2t^2(u')^2 + 2t^2uu'' + 6t uu'.
\]
which upon substituting for $u, u'$ and $u''$ become 
\[
M_1(t) = 2t^3(t-1)^2 \ti M_1(t) /(t \log t )^6, \; \text{and} \; M_2(t) = \ti M_2(t)/t^{2}(\log t)^{4}
\]
where 
\begin{align*}
\ti M_1(t) &= \left((t-1)^2 - t (\log t)^2  \right), \; \text{and} \\
\ti M_2(t) &= (t-1)^2(\log t)^2 + 2 (\log t)^2 + 6(t-1)^2 - 8(t-1) \log t + 2(t-1)(\log t)^2 - 4(t-1)^2 \log t 
\end{align*}
It can now be seen that $M_2(t)/(tu^2)^3 \ra 0$ as $t \ra 0$ and $M_1(t)/(tu^2)^3 \ra 2$ and hence $f(t) \ra -4$ as $t \ra 0$.

\medskip

Further, both $\ti M_1$ and $\ti M_2$ approach $+\infty$ as $t \ra \infty$ and hence $M_1, M_2 > 0$ for large $t$. Since $u > 0, u' \le 0$ and $u'' \ge 0$ on $(0, \infty)$, each of the terms $6tuu', 2t^2u^3u'$ and $-2t^3(uu')^2$  are non-positive, and hence for large $t$,
\[
0 \le M_2(t) = u^2 + 2t^2(u')^2 + 2t^2 uu'' + 6tuu' \le u^2 + 2t^2(u')^2 + 2t^2 uu'' 
\]
and 
\[
0 \le M_1(t) = 2t^3u^3u'' + 2t^2u^3u' - 2t^3(uu')^2 \le 2t^3u^3u''. 
\]
Both $(u^2 + 2t^2(u')^2 + 2t^2 uu'')/(tu^2)^3$ and $2t^3u^3u''/(tu^2)^3 \ra 0$ as $t \ra 0$ and thus, $f \ra 0$ as $t \ra \infty$. 

\medskip
 
Finally, the assertion that $f \le 0$ on $(0, \infty)$ is equivalent to the subharmonicity of $\log h$ on $\mbb C^{\ast}$. For this, it suffices to show that $\pa \ov \pa \log (\vert z \vert^2 u^2(\vert z \vert^2)  = \pa \ov \pa \log u^2(\vert z \vert^2) \ge 0$ on $\mbb C^{\ast}$. Now compute
\begin{align*}
\pa \ov \pa \log u(\vert z^2 \vert)  & = u^{-2} \left( uu' + \vert z \vert^2 \left( uu'' - (u')^2 \right) \right) \\
                                                         & = t^{-3} (\log t)^{-4} \left( (t-1)^2 - t (\log t)^2 \right). 
\end{align*}
Ignoring the positive multiplicative factor $t^{-3}$, the sign of $\pa \ov \pa  \log u(\vert z^2 \vert)$ is determined by 
\[
(\log t)^{-4} \left( (t-1)^2 - t (\log t)^2 \right).
\]
Writing $t = e^x$, where $x \in \mbb R$ now, reduces this to
\[
x^{-4}(e^x - 1)^2 - x^{-2}e^x
\]
whose value near $x = 0$ is positive since its limit at $x = 0$ is $1/12$. Away from $x = 0$, the positivity of this expression is equivalent to the assertion that $\sinh^2(x/2) \ge (x/2)^2$ which is evidently true. 

\medskip

Thus, the Grauert metric on $\mbb C^{\ast}$ is an example of a non-positively curved conformal metric whose limiting behaviour is described above. Given that $\kappa \ra -4$ near $z = 0$, the fact that $\kappa \ra 0$ for large $\vert z \vert$ is not surprising and in fact, is expected. Indeed, if $\kappa$ were to approach a negative constant as $\vert z \vert \ra \infty$, then $g$ would be an example of a metric on $\mbb C^{\ast}$ that is non-positively curved everywhere and, in addition, satisfies $\kappa \le -a < 0$ outside a compact set in $\mbb C^{\ast}$. This would contradict a version of the Schwarz lemma -- see Troyanov \cite{Tr}.

\medskip

\no (ii) When $n \ge 2$, working directly with the curvature tensor to understand the limiting behaviour near the origin and for large values of $\vert z \vert$ involves the inversion of the matrix 
$((g_{\al \ov \be}))$. To circumvent this, the approach of working with the foliation induced by a suitable holomorphic vector field has been adopted. The question of which holomorphic vector field to work with then arises, and as the calculations show, the only condition to keep in mind is that $\mathcal S$ be non-empty. The choice of $X$ as in (1.1), which is a linear saddle-node example, clearly satisfies this condition. In this case, the coordinate axes are separatrices and Theorem 1.1  shows that the leaf-wise curvature $\kappa \ra +\infty$ along $\mathcal S$ near the origin. All these considerations are related to the paradigm of identifying Hermitian metrics which are negatively curved when restricted to the leaves of a given foliation -- for example, see \cite{CNS}, \cite{Ca}, \cite{CM}, \cite{Gl}, \cite{Ne1}, \cite{Ne2}, \cite{Shc} and \cite{Ve}.


\end{document}